\documentclass[graybox,envcountsame,envcountsect]{svmult}
\usepackage{amsmath,amscd,amssymb}
\usepackage{latexsym}
\usepackage{graphicx}
\usepackage{enumerate}

\usepackage{hyperref}

\newcommand{\C}{{\mathbb{C}}}
\newcommand{\bbC}{{\mathbb{C}}}

\newcommand{\DD}{{\mathbb{D}}}

\newcommand{\R}{{\mathbb{R}}}
\newcommand{\bbR}{{\mathbb{R}}}
\newcommand{\Z}{{\mathbb{Z}}}

\newcommand{\al}{\alpha}
\newcommand{\be}{\beta}

\newcommand{\ga}{\gamma}
\newcommand{\la}{\lambda}
\newcommand{\Om}{\Omega}

\newcommand{\fre}{{\mathfrak{e}}}
\newcommand{\frf}{{\mathfrak{f}}}

\newcommand{\bi}{\bibitem}

\newcommand{\lb}{\label}

\newcommand{\Oh}{O}

\newcommand{\st}{:}

\newcommand{\pd}{\partial}

\newcommand{\cvh}{\text{\rm{cvh}}}
\newcommand{\inte}{\text{\rm{int}}}

\DeclareMathOperator{\ca}{\mathrm{C}}

\spnewtheorem*{remarks}{Remarks}{\bf}{}
\spnewtheorem{openproblem}[theorem]{Open Problem}{\bf}{}

\spnewtheorem{conjecture}[theorem]{Conjecture}{\bf}{}

\begin{document}

\title*{Widom Factors and Szeg\H{o}--Widom Asymptotics, a Review}
\author{Jacob S.~Christiansen, Barry Simon, and Maxim~Zinchenko}

\institute{
Jacob S.~Christiansen \at
Centre for Mathematical Sciences, Lund University, Box 118, 22100 Lund, Sweden;
\email{stordal@maths.lth.se}\\
Research supported by VR grant 2018-03500 from the Swedish Research Council and in part by DFF research project 1026-00267B from the Independent Research Fund Denmark.
\and
Barry Simon \at Departments of Mathematics and Physics, Mathematics 253-37, California Institute of Technology, Pasadena, CA 91125, USA;
\email{bsimon@caltech.edu}\\
Research supported in part by Israeli BSF Grant No.~2020027.
\and
Maxim~Zinchenko \at
Department of Mathematics and Statistics, University of New Mexico,
Albuquerque, NM 87131, USA;
\email{maxim@math.unm.edu}\\
Research supported in part by Simons Foundation grant CGM-581256.
}

\maketitle
\abstract*{}
\abstract{
We survey results on Chebyshev polynomials centered around the work of H.~Widom. In particular, we discuss asymptotics of the polynomials and their norms and general upper and lower bounds for the norms. Several open problems are also presented.
}

\bigskip


\medskip



\vspace{3em}
\emph{Dedicated with great respect to the memory of Harold Widom, 1932--2021.}
\vspace{2em}

\numberwithin{equation}{section} 

\section{Introduction} \lb{s1}

Let $\fre\subset\C$ be a compact, not finite set and denote by $$\Vert f \Vert_\fre:=\sup_{z\in\fre} \vert f(z) \vert$$ the supremum norm of a continuous, complex-valued function $f$ on $\fre$.  A classical problem in approximation theory is, for every $n\geq 1$, to find the unique monic degree $n$ polynomial, $T_n$, which minimizes $\Vert P \Vert_\fre$ among all monic degree $n$ polynomials, $P$.  The resulting sequence is called the \emph{Chebyshev polynomials} of $\fre$.

By the maximum principle, we may assume that $\fre$ is polynomially convex.  This means that $\Omega:=(\C\cup\{\infty\})\setminus\fre$ is connected so that $\fre$ has no inner boundary.

It is only in the case of $\fre$ being a (possibly elliptical) disk or a line segment that explicit formulas for all $T_n$'s are available.  The Chebyshev polynomials of the unit disk are simply $T_n(z)=z^n$, while the ones for the interval $[-1, 1]$ (or any ellipse with foci at $\pm1$) are given by $$T_n(x)=2^{-n+1}\cos(n\theta),$$ where $x=\cos\theta$.

In addition to this, there are certain sets generated by polynomials (such as lemniscates and Julia sets) for which a subsequence of $T_n$ can be written in closed form.
For general $\fre$, however, the best one can hope for is to determine the asymptotic behavior of $T_n$.  In this article we seek to present what is known about the asymptotics of Chebyshev polynomials.  Had it not been for Widom's landmark paper \cite{Wid69}, there probably wouldn't be much to say.

To get started, we briefly introduce some notions from potential theory (see, e.g., \cite{ArmGar01,Hel09,Lan72,MF06,Ran95} for more details).  Let $\ca(\fre)$ denote the logarithmic capacity of $\fre$.  When $\fre$ is non-polar (i.e., $\ca(\fre)>0$), we denote by $d\rho_\fre$ the equilibrium measure of $\fre$ and by $G:=G_\fre$ the Green's function of $\fre$.  These are closely linked by the relation
\begin{equation}
 G(z)=-\log \bigl[ \ca(\fre) \bigr]+\int\log\vert z-x \vert d\rho_\fre(x).
\end{equation}
For subsets $\frf\subset\fre$, we shall also refer to $\rho_\fre(\frf)$ as the harmonic measure of $\frf$. 
The set $\fre$ is called {regular} if $G$ vanishes at all points of $\fre$ (equivalently, $G$ is continuous on all of $\C$).

The general results for Chebyshev polynomials are few, but important. Szeg\H{o} \cite{Sze1924} showed that
\begin{equation}
\label{Sze1}
 \Vert T_n \Vert_\fre \geq \ca(\fre)^n.
\end{equation}
This applies to all compact sets $\fre\subset\C$ and is optimal since equality occurs for all $n$ when $\fre$ is a disk.  When $\fre\subset\R$, Schiefermayr \cite{Sch08} improved upon \eqref{Sze1} by showing that
\begin{equation}
\label{Sch}
 \Vert T_n \Vert_\fre \geq 2\ca(\fre)^n, \quad n\geq 1,
\end{equation}
which is again optimal (take $\fre$ to be an interval).  Szeg\H{o} \cite{Sze1924}, using prior results of Faber \cite{Fab1919} and Fekete \cite{Fek1923}, also proved the following asymptotic result:
\begin{equation}
\label{Sze2}
 \lim_{n\to\infty} \Vert T_n \Vert_\fre^{1/n} = \ca(\fre).
\end{equation}
This certainly puts a growth restriction on $\Vert T_n \Vert_\fre$ but is not strong enough to force $\Vert T_n \Vert_\fre/\ca(\fre)^n$ to be bounded.  We shall discuss which extra assumptions on $\fre$ may imply this in Sections~\ref{s2}--\ref{s3}.

The polynomials themselves also obey $n$th root asymptotics.  For a non-polar compact set $\fre\subset\C$, we have that
\begin{equation}
\label{root}
 \vert T_n(z) \vert^{1/n} \to \ca(\fre)\exp \bigl[ G(z) \bigr]
\end{equation}
uniformly on any closed set disjoint from $\cvh(\fre)$, the convex hull of $\fre$.  This result is implicitly in Widom \cite{Wid67}, where he shows that all zeros of $T_n$ must lie in $\cvh(\fre)$ before proceeding to the asymptotics.  See also Ullman \cite{Ull59} and Saff--Totik \cite[Chap.~III]{SafTot97}.

``All asymptotic formulas have refinements,'' quoting the introduction of \cite{Wid69}.  And this is precisely what we aim at, just as Widom did.  As \eqref{Sze1}--\eqref{root} suggest, it is natural to scale $T_n$ by a factor of $\ca(\fre)^n$.  We shall study the limiting behavior of the so-called \emph{Widom factors}
\begin{equation}
\label{Wf}
 W_n(\fre):=\Vert T_n \Vert_\fre / \ca(\fre)^n.
\end{equation}
If this scaled version of the norms does not have a limit, can we then at least single out the possible limit points?
Regarding the polynomials $T_n$, we aim at strong asymptotics or what we shall refer to as \emph{Szeg\H{o}--Widom asymptotics}.

The first result in this direction goes back to Faber \cite{Fab1919}.  When $\fre$ is a closed Jordan region, there is a Riemann map of $\Om$ onto the unit disk, $\DD$.  We uniquely fix this map, $B$, by requiring that
\begin{equation}
\label{B}
 B(z)=\frac{\ca(\fre)}{z}+\Oh(1/z^2)
\end{equation}
near $\infty$.  Assuming that $\partial\fre$ is analytic, Faber showed that $W_n(\fre)\to 1$ and, more importantly, that
\begin{equation}
\label{Faber}
 \frac{T_n(z) B(z)^n}{\ca(\fre)^n} \to 1
\end{equation}
uniformly for $z$ in a neighborhood of $\overline{\Om}$.

The picture changes completely when $\fre$ consists of more than one component.  In his work on Chebyshev polynomials of two intervals, Akhiezer \cite{Ach1, Ach2} proved that either $W_n(\fre)$ is asymptotically periodic or else the set of limit points of $W_n(\fre)$ fills up an entire interval.  But it was only Widom \cite{Wid69} who lifted the theory to $\fre$ being a union of disjoint compact subsets of $\C$ and developed a framework to distinguish between periodicity and almost periodicity.

In replacement of the Riemann map, we introduce (on $\Om$) a multivalued analytic function $B:=B_\fre$ which is determined by
\begin{equation}
\label{Be}
 \vert B(z) \vert = \exp\bigl[ -G(z) \bigr]
\end{equation}
and \eqref{B} near $\infty$.  One can construct this $B$ using the fact that $-G$ is locally the real part of an analytic function whose exponential ($=B$) can be continued along any curve in $\Om$.  By the monodromy theorem, the continuation is the same for homotopic curves and, due to \eqref{Be}, going around a closed curve $\ga$ can only change $B$ by a phase factor.  Hence there is a character $\chi_\fre$ of the fundamental group $\pi_1(\Om)$ so that going around $\ga$ changes $B$ by $\chi_\fre([\ga])$.  More explicitly, if $\ga$ winds around a subset $\frf\subset\fre$ and around no other points of $\fre$, then the multiplicative change of phase of $B$ around $\ga$ is given by
\begin{equation}
 \exp\bigl[-2\pi i\rho_\fre(\frf)\bigr].
\end{equation}

In line with Faber and \eqref{Faber}, Widom looked at ${T_n(z) B(z)^n}/{\ca(\fre)^n}$ for the ``new'' $B$ and noted that its character $\chi^n_\fre$ only has a limit when $\chi_\fre$ is trivial (i.e., $\Om$ is simply connected).  So there is no hope of finding a pointwise limit except when $\fre$ just has one component.  Widom's  stroke of genius was to find a good candidate for the asymptotics when $\fre$ has several components.  For every character $\chi$ in $\pi_1(\Om)^*$ there exists a so-called \emph{Widom minimizer} which we shall denote by $F_\chi$.  This is the unique element of $H^\infty(\Om, \chi)$ (i.e., the set of bounded analytic $\chi$-automorphic function on $\Om$) with $F_\chi(\infty)=1$ and for which
\begin{equation}
\label{F}
 \Vert F_\chi \Vert_\infty = \inf\bigl\{ \Vert h \Vert_\infty \st h\in H^\infty(\Om,\chi), \, h(\infty)=1 \bigr\}.
\end{equation}
Writing $F_n$ as shorthand notation for $F_{\chi^n}$, the \emph{Widom surmise} is the notion that
\begin{equation}
\label{W-surmise}
 \frac{T_n(z) B(z)^n}{\ca(\fre)^n}-F_n(z)\to 0.
\end{equation}
When it holds uniformly on compact subsets of the universal cover of $\Om$, we say that $\fre$ has {Szeg\H{o}--Widom asymptotics}.

Widom \cite{Wid69} proved that one has this type of asymptotics when $\fre$ is a finite union of disjoint Jordan regions with smooth boundaries and conjectured that this should also hold for finite gap sets (in $\R$).  A main result of \cite{CSZ1} was to settle this conjecture.  By streamlining the method of proof, this was then extended to a large class of infinite gap sets in \cite{CSYZ2} (see Section~\ref{s2} for further details).

The framework of characters is also useful when describing the fluctuation of $W_n(\fre)$.  In \cite{Wid69}, Widom proved that
\begin{equation}
\label{W-asympt}
 W_n(\fre) / \Vert F_n \Vert_\infty \to 1
\end{equation}
for finite unions of disjoint Jordan regions and established the counterpart (with $1$ replaced by $2$ on the right-hand side) for finite gap sets.  The behavior of $\Vert F_n \Vert_\infty$ very much depends on the character $\chi_\fre$.  If $\chi_\fre^n=1$ for some $n$, then the sequence is periodic (with period at most $n$) and otherwise it is merely \emph{almost periodic}.  This is precisely the pattern that Akhiezer discovered for two intervals.  We shall discuss the possible limit points in more detail in Section~\ref{s2}.

The paper is organized as follows.  In Section~\ref{s2} we discuss bounds and asymptotics for Chebyshev polynomials of compact subsets of the real line.  Then in Section~\ref{s3} we survey similar results for Chebyshev and weighted Chebyshev polynomials of subsets of the complex plane, including results on the asymptotic distribution of zeros.  Open problems are formulated along the way.

We would be remiss if not mentioning related problems, such as the Ahlfors problem \cite{EicYud18}, and similar classes of polynomials or functions, for instance, residual polynomials \cite{CSZ5,Yud99} and rational Chebyshev functions \cite{EMY22}.  But to consider the subject in more depth, we decided to merely focus on the Chebyshev problem.


\section{Real Chebyshev Polynomials} \lb{s2}

As we shall see, there is a rather complete theory for Chebyshev polynomials of compact sets $\fre\subset\R$.  This is in part due to what is called \emph{Chebyshev alternation}.  We say that $P_n$, a real degree $n$ polynomial, has an alternating set in $\fre$ if there exists $n+1$ points in $\fre$, say $x_0<x_1<\ldots<x_n$, so that
\begin{equation}
 P_n(x_j) = (-1)^{n-j} \Vert P_n \Vert_\fre.
\end{equation}
The alternation theorem gives the following characterization of the $n$th Chebyshev polynomial of $\fre$:
\emph{$T_n$ always has an alternating set in $\fre$ and, conversely, any monic degree $n$ polynomial with an alternating set in $\fre$ must be equal to $T_n$}.

This result, in turn, has consequences for the zeros of $T_n$.  Not only do all of them lie in $\cvh(\fre)$, but any gap of $\fre$ (i.e., a bounded component of $\R\setminus\fre$) contains at most one zero of $T_n$.  The alternating set need not be unique and usually isn't.  However, it always contains the endpoints of $\cvh(\fre)$.  See, e.g., \cite{CSZ1} for proofs and more details.

We now turn the attention to the Widom factors which were introduced in \eqref{Wf}.  By \cite{Sch08} we always have $W_n(\fre)\geq 2$ and, as proven in \cite{CSZ3}, equality occurs for $n=km$ (with $m\geq 1$) precisely when
\begin{equation}
\label{preimg}
 \fre=P^{-1}\bigl( [-2,2] \bigr)
\end{equation}
for some degree $k$ polynomial, $P(z)=cz^k+\mbox{lower order terms}$.  In that case we actually have $T_{km}=(P/c)^m$ for all $m\geq 1$, an observation that essentially goes back to Faber \cite{Fab1919}.  It also follows that equality holds for all $n$ if and only if $\fre$ is an interval.  A stronger and related result of Totik \cite{Tot14} states that if $\lim_{n\to\infty} W_n(\fre)=2$, then $\fre$ must be an interval.

Interestingly, the sets that appear in \eqref{preimg} are not only of interest for the lower bound; they play a key role in the theory.  For $\fre\subset\R$, we introduce the so-called \emph{period-$n$ sets}, $\fre_n$, (aka \emph{$n$-regular sets} \cite{SY}) by
\begin{equation}
\label{En}
 \fre_n:=T_n^{-1}\bigl( [-\Vert T_n \Vert_\fre, \Vert T_n \Vert_\fre] \bigr). 
\end{equation}
Clearly, $T_n$ is also the Chebyshev polynomial of $\fre_n\supset\fre$ and furthermore we have that
\begin{equation}
\label{equal}
\Vert T_n \Vert_\fre = 2\ca(\fre_n)^n.
\end{equation}
Due to alternation we can write any period-$n$ set as
\begin{equation}
\label{ab}
 \fre_n=\bigcup_{j=1}^n \, [\al_j, \be_j],
\end{equation}
where $\al_1<\be_1\leq \ldots \leq \al_n<\be_n$ are the solutions of $T_n(x)=\pm \Vert T_n \Vert_\fre$.  So $T_n$ is strictly monotone on each of the bands $[\al_j, \be_j]$ and $\fre_n\subset\cvh(\fre)$.  Note that $\al_1$ and $\be_n$ always belong to $\fre$ while for $j=1,\ldots,n-1$, at least one of $\be_j$ and $\al_{j+1}$ must lie in $\fre$.  Therefore, any gap of $\fre$ can at most overlap with one of the bands of $\fre_n$.

The period-$n$ sets are well suited for potential theory.  For instance, the Green's function and equilibrium measure of $\fre_n$ are explicitly given by
\begin{equation}
 G_n(z)=\frac{1}{n}\log\Biggl\vert \frac{\Delta_n(z)}{2}+\sqrt{\biggl(\frac{\Delta_n(z)}{2}\biggr)^2-1\,} \Biggr\vert
\end{equation}
and
\begin{equation}
 d\rho_n(x)=\frac{1}{\pi n} \frac{\vert \Delta_n'(x)\vert}{\sqrt{4-\Delta_n(x)^2}} \,dx, \quad x\in\fre_n,
\end{equation}
where $\Delta_n$ is defined by
\begin{equation}
\label{Del n}
 \Delta_n(z):=2 T_n(z)/\Vert T_n \Vert_\fre.
\end{equation}
In particular, each band of $\fre_n$ has $\rho_n$-measure 
$1/n$.  See, e.g., \cite{CSZ1} for proofs and further details.

Comparing the Green's functions for $\fre$ and $\fre_n$ at $\infty$, and letting $\{ K_j \}$ account for the gaps of $\fre$, we see that
\begin{equation}
 \log\bigl[ \ca(\fre_n)/\ca(\fre) \bigr] = \int_{\fre_n} \bigl[G(x)-G_n(x)\bigr] d\rho_n(x) \leq \frac{1}{n}\sum_j \max_{x\in K_j} G(x)
\end{equation}
which combined with \eqref{equal} then yields
\begin{equation}
 \Vert T_n \Vert_\fre/ \ca(\fre)^n \leq 2 \sum_j \max_{x\in K_j} G(x).
\end{equation}
This observation leads to an upper bound on $W_n(\fre)$ for a large class of compacts sets $\fre\subset\R$.  When $\fre$ is regular (for potential theory), the Green's function vanishes at all endpoints of the $K_j$'s and since $G$ is also concave on the gaps, it attains its maximum on $K_j$ at the unique critical point, $c_j$, in that gap.  A regular compact set $\fre\subset\R$ (or $\C$) is called a \emph{Parreau--Widom set} (in short, PW) if
\begin{equation}
\label{PW}
 PW(\fre):=\sum_j G(c_j)<\infty,
\end{equation}
where the sum is over all points $c_j\in\R\setminus\fre$ for which $\nabla G(c_j)=0$.  Such sets are known to have positive Lebesgue measure (see, e.g., \cite{JSC} for details).  One of the main results of \cite{CSZ1} that we have now deduced is the following:
\begin{theorem}
\label{TW}
If $\fre\subset\R$ is a PW set, then the Widom factors are bounded.  Explicitly, we have that
\begin{equation}
\label{TWub}
 \Vert T_n \Vert_\fre \leq 2\exp\bigl[PW(\fre)\bigr] \ca(\fre)^n.
\end{equation}
\end{theorem}
\begin{remarks}
$(i)$ Sets which obey \eqref{PW} were introduced by Parreau \cite{Par58} in the context of Riemann surfaces.  They later appeared in Widom's work on multi-valued analytic functions \cite{Wid71a, Wid71b} and the name was coined by Hasumi in his monograph \cite{Hasu}.\\
$(ii)$ Examples of PW sets include finite gap sets but also sets that are homogeneous in the sense of Carleson \cite{Car81}, e.g., fat Cantor sets.\\  $(iii)$ Upper bounds of the form $\Vert T_n \Vert_\fre \leq K\cdot\ca(\fre)^n$ are also referred to as \emph{Totik--Widom} bounds. Here $K>0$ is a constant that does not depend of $n$.
\end{remarks}

As alluded to in the introduction, the Widom factors are not always bounded.  It was proven in \cite{BGH83} that they are unbounded when $\fre$ is the Julia set of $(z-\la)^2$ and $\la>2$.  Interestingly, $W_{2^n}(\fre)$ is bounded (in fact, constant) in that case, while $W_{2^n-1}(\fre)\to\infty$.  There are more elaborate examples of very thin Cantor-type sets for which $W_n(\fre)$ grows subexponentially of any order, see Goncharov--Hatino\v{g}lu \cite{GH15} for details.  But it is not known if the Widom factors of, e.g., the middle $3$rd Cantor set are bounded.  The best result in this direction, due to Andrievskii \cite{And17}, states that when $\fre\subset\R$ is uniformly perfect there exists a constant $c>0$ such that $W_n(\fre)=\Oh(n^c)$.  We pose the following question:
\begin{openproblem}
Does there exist a Lebesgue measure zero set or merely a non-PW set $\fre\subset\R$ for which the Widom factors are bounded?
\end{openproblem}

It remains to consider the fluctuation and possible limit points of $W_n(\fre)$.  We shall do so in conjunction with the asymptotics of the polynomials.

Let us start by explaining, following \cite{CSYZ2, CSZ1}, how one can establish Szeg\H{o}--Widom asymptotics.  Since every band of $\fre_n$ has $\rho_n$-measure $1/n$, the $n$th power of $B_n:=B_{\fre_n}$ 
is single-valued.  In fact,
\begin{equation}
 B_n(z)^{\pm n}=\frac{\Delta_n(z)}{2}\mp \sqrt{\biggl(\frac{\Delta_n(z)}{2}\biggr)^2-1\,}
\end{equation}
with $\Delta_n$ as in \eqref{Del n}.  It follows that
\begin{equation}
\label{key}
 \frac{2T_n(z)}{\Vert T_n \Vert_\fre}=B_n(z)^n + B_n(z)^{-n}
\end{equation}
and this is the key formula we need.  As a side remark we note that when $\fre=[-1,1]$, \eqref{key} corresponds to the familiar formula
\begin{equation}
 T_n(z)=\frac{1}{2^n} \biggl( \Bigl( z-\sqrt{z^2-1} \Bigr)^n + \Bigl( z+\sqrt{z^2-1} \Bigr)^n \biggr).
\end{equation}
The idea is now to recast \eqref{key} in the form
\begin{equation}
\label{idea}
 \frac{T_n(z) B(z)^n}{\ca(\fre)^n}=\Bigl( 1+ B_n(z)^{2n} \Bigr) \frac{M_n(z)}{M_n(\infty)},
\end{equation}
where
\begin{equation}
 M_n(z)=B(z)^n / B_n(z)^n.
\end{equation}
Since $\sup_{n, \, z\in K} \vert B_n(z) \vert <1$ on any compact subset $K$ of the universal cover of $\Om$, the task is reduced to proving that
\begin{equation}
 F_n(z)-M_n(z)/M_n(\infty) \to 0
\end{equation}
and this can be done by controlling the limit points of $M_n$.

In order to go beyond finite gap sets, some issues have to be sorted out.  First of all, for which infinite gap sets do the Widom minimizers at all exist and are they unique?  Secondly, how are the limit points of $M_n$ related to the Widom minimizers and do these minimizers depend continuously on the character so that one can pass to the limit along convergent subsequences?

The answer to both of the above questions are rooted in Widom's work.  In \cite{Wid71b}, he proved that \eqref{PW} holds if and only if there is a nonzero element in $H^\infty(\Om,\chi)$ for every $\chi\in\pi_1(\Om)^*$.  Hence, by compactness, Widom minimizers exist for all PW sets.  Uniqueness requires a separate argument for which we refer the reader to \cite{CSYZ2} and \cite{VY14}.  Note also that Theorem \ref{TW} implies $\vert M_n \vert \leq 1$ in the PW regime.  So limit points do exist in that setting by Montel's theorem.

To proceed with the analysis, it is instructive to also consider the problem dual to \eqref{F}.  The function $Q_\chi\in H^\infty(\Om, \chi)$ which satisfies
\begin{equation}
 Q_\chi(\infty) = \sup\bigl\{ g(\infty) \st g\in H^\infty(\Om,\chi), \, \Vert g \Vert_\infty=1, \, g(\infty)>0 \bigr\}
\end{equation}
is called the \emph{dual Widom maximizer}.  Clearly, we have
\begin{equation}
 Q_\chi=F_\chi/ \Vert F_\chi \Vert_\infty, \quad
 F_\chi=Q_\chi/ Q_\chi(\infty), \quad
 Q_\chi(\infty)=1/ \Vert F_\chi \Vert_\infty
\end{equation}
and therefore the two problems either both or neither have unique solutions.  By controlling the zeros of $T_n$ in gaps of $\fre$ and using the fact that
\begin{equation}
 \vert M_n(z) \vert=\exp\biggl[ -n \int_{\cup_j K_j} G(x,z) d\rho_n(x) \biggr],
\end{equation}
one can prove that the limit points of $M_n$ are dual Widom maximizers.  More precisely, the approach of \cite{CSYZ2} reveals that limit points of $M_n$ are Blaschke products with at most one zero per gap of $\fre$ and such character automorphic products are indeed dual Widom maximizers.

As for the final issue, Widom \cite{Wid71a} noted that ``It is natural to ask (and important to know) whether $Q_\chi$ is continuous as a function of $\chi$ on the compact group $\pi_1(\Om)^*$.''  He pointed out that this can easily fail to hold (e.g., if $\fre$ has isolated points) but was not able to characterize those sets for which we have continuity.  Years later, this was settled by Hayashi and Hasumi (see \cite{Hasu, Hay}).  Continuity in $\chi$ is equivalent to having a so-called \emph{direct Cauchy theorem} (DCT) on $\Om$.  There seems to be no obvious geometric interpretation of this DCT property;  while it may fail for a general PW set, it always holds when $\fre$ is homogeneous (see, e.g., \cite{Yud2011} for further details).

We should point out that DCT is responsible for the almost periodic behavior of the Widom minimizers.  That is,
\begin{equation}
\label{ap1}
 n\mapsto \Vert F_n \Vert_\infty \mbox{ is an almost periodic function}
\end{equation}
and
\begin{align}
\label{ap2}
 \notag
 n\mapsto F_n(z) & \mbox{ is almost periodic uniformly for $z$ in compact} \\
 & \mbox{ subsets of the universal cover of $\Om$.}
\end{align}
Recall namely that $n\mapsto x_n$ is almost periodic precisely when $\{x_n\}$ is the orbit of a continuous function on a torus (possibly of infinite dimension).  Since the character group $\pi_1(\Om)^*$ is topologically a torus, we are led to \eqref{ap1}--\eqref{ap2}.

After this extended discussion, we are now ready to formulate the main result of \cite{CSYZ2}.
\begin{theorem}
\label{SW as}
If $\fre\subset\R$ is a PW set and obeys the DCT condition, then the Chebyshev polynomials of $\fre$ have strong Szeg\H{o}--Widom asymptotics.
That is, the Widom surmise \eqref{W-surmise} and \eqref{ap1}--\eqref{ap2} hold.  Moreover,
\begin{equation}
\label{norm a}
 \lim_{n\to\infty} \frac{\Vert T_n \Vert_\fre}{\ca(\fre)^n \Vert F_n \Vert_\infty} = 2.
\end{equation}
\end{theorem}

\begin{remarks}
$(i)$  The additional word ``strong'' is used here to include the almost periodicity of \eqref{ap1}--\eqref{ap2}.\\
$(ii)$  The last statement also follows from \eqref{idea} by noting that $\sup_{z\in\Om} \vert 1+B_n(z)^{2n} \vert=2$ since there are points $x\in\fre_n$ with $\vert B_n(x) \vert=1$.
\end{remarks}

The above theorem enables us to shed more light on the fluctuation of the Widom factors.  {When $\fre\subset\R$ is a PW set with DCT, the function $n\mapsto W_n(\fre):=\Vert T_n \Vert_\fre / \ca(\fre)^n$ is asymptotically almost periodic}.  The set of limit points may or may not fill up the entire interval between the lower bound (=$2$) and the upper bound from Theorem \ref{TW}.  Generically, it will (as explained in \cite{CSZ3}; see also below) but this is not the case when, for instance, $\fre$ is a period-$n$ set.  For in that case, we have $\chi_\fre^n=1$ and the function in \eqref{ap1} becomes periodic.

Following \cite{CSYZ2}, we say that $\fre$ has a \emph{canonical generator} if the orbit $\{ \chi_\fre^n \}_{n\in\Z}$ is dense in $\pi_1(\Om)^*$.  This holds if and only if for all decompositions $\fre=\fre_1\cup\ldots\cup\fre_l$ into disjoint closed sets and rational numbers $\{ q_j \}_{j=1}^{l-1}$  (not all zero), we have that
\begin{equation}
 \sum_{j=1}^{l-1} q_j \rho_\fre(\fre_j) \neq 0 \; \mbox{ (mod $1$)}.
\end{equation}
In particular, a finite gap set has a canonical generator precisely when the harmonic measures of the bands are rationally independent (except that they sum to $1$).  One can show that the property of having a canonical generator is generic (see \cite{CSYZ2} for details) and it implies that any number $\geq 2$ and $\leq 2\exp\bigl[ PW(\fre) \bigr]$ is a limit point of $W_n(\fre)$.

To end this section, we return to the open problem formulated a few pages ago.  It was proven in \cite{CSYZ2} that if $\fre$ has a canonical generator and obeys a Totik--Widom bound (as in Theorem \ref{TW}), then it must be a PW set.  This provides some evidence that the answer could be in the negative.
However, there are also results pulling in the opposite direction.  While $\liminf W_n(\fre)=2$ for any PW set with DCT (as proven in \cite{CSZ5}), it is not always the case that $\limsup$ is equal to $2\exp\bigl[PW(\fre)\bigr]$ when $\fre$ is a period-$n$ set and $n\geq 2$.  For instance, one can prove that strict inequality applies when $\fre$ is the degenerate period-$3$ set $$\bigl[-\sqrt{3}, 0\bigr] \cup \bigl[\sqrt{3}, 2\bigr].$$  We thus wonder if some cleverly arranged limit of period-$n$ sets could provide an example of a non-PW set with bounded Widom factors.

\section{Complex Chebyshev Polynomials} \lb{s3}

In this section we consider Chebyshev and weighted Chebyshev polynomials for compact subsets of the complex plane. In particular, we discuss Widom's contribution to the subject as well as several recent refinements.

Throughout the section we will assume that $\ca(\fre)>0$
and let $w$ be a nonnegative upper semi-continuous weight function on $\fre$ (this ensures that $w$ is bounded) which is nonzero at infinitely many points of $\fre$.  Under these assumptions, there exists for each $n\ge1$ a unique weighted Chebyshev polynomial $T_{n,w}:=T_{n,w}^{(\fre)}$ that minimizes $\|w T_{n,w}\|_\fre$ among monic polynomials of degree $n$.

In \cite{Wid67}, Widom proved that one has root asymptotics analogous to \eqref{Sze2} for a fairly general class of monic extremal polynomials which, in particular, includes the $L^p(wd\rho_\fre)$-extremal polynomials for $0<p<\infty$ and weights $w$ satisfying $w>0$ $d\rho_\fre$-a.e.
This type of asymptotics for the $L^1(wd\rho_\fre)$-extremal polynomials, $P_n$, combined with \eqref{Sze2} and the two-sided estimate 
\begin{align}
  \|P_n\|_{L^1(wd\rho_\fre)} \le \|T_{n,w}\|_{L^1(wd\rho_\fre)} \le \|w T_{n,w}\|_\fre \le \|w T_n\|_\fre \le \|w\|_\fre\|T_n\|_\fre
\end{align}
yields that the weighted Chebyshev polynomials obey the root asymptotics 
\begin{align}
  \lim_{n\to\infty} \|w T_{n,w}\|_\fre^{1/n} = \ca(\fre)
\end{align}
whenever $w>0$ $d\rho_\fre$-a.e.

As explained below, there is also a lower bound, an asymptotic upper bound, and strong asymptotics for the weighted Chebyshev polynomials under an additional assumption on the weight function $w$, namely the so-called Szeg\H{o} condition
\begin{align} \lb{Sz}
  S(w) = \exp\left[\int\log w(z)\,d\rho_\fre(z)\right] > 0.
\end{align}
A generalization of Szeg\H{o}'s lower bound \eqref{Sze1} to the weighted case was observed for finite unions of Jordan regions by Widom \cite[Sect.~8]{Wid69} and extended to general non-polar compact sets $\fre\subset\bbC$ in \cite{NSZ21}.  It relies on \eqref{Sz} and states that
\begin{align}
\label{wLB}
  \|w T_{n,w}\|_\fre \ge S(w)\ca(\fre)^n.
\end{align}
In addition, it was shown in \cite{NSZ21} that unlike the unweighted case, this lower bound is sharp even for real sets $\fre$ (cf. \eqref{Sch}). Moreover, equality in \eqref{wLB} occurs for some $n$ if and only if there exists a monic polynomial, $P_n$, of degree $n$ such that $P_n(z)=0$ implies $G(z)=0$ and $w(z)|P_n(z)|=\|w P_n\|_\fre$ for $d\rho_\fre$-a.e.\ $z\in\fre$, in which case $T_{n,w}=P_n$.


Next, we turn to upper bounds.  A collection of very general bounds for unweighted Chebyshev polynomials were obtained by Andrievskii \cite{And16,And17} and Andrievskii--Nazarov \cite{AN19}.  See also Totik--Varga \cite{ToVa15}.  In particular, it was shown that if $\fre\subset\bbC$ is a finite union of  quasiconformal arcs and/or Jordan regions bounded by quasiconformal curves (aka quasidisks), then a Totik--Widom upper bound $$\|T_n\|_\fre \le K\cdot \ca(E)^n$$ holds with some constant $K$. This result includes a large class of regions with pathological boundaries, for example, the Koch snowflake.  In addition, in the absence of any smoothness it was shown that for compact sets $\fre$ with finitely many components, the Widom factors $W_n(\fre)$ can grow at most logarithmically in $n$. Yet, in this setting no example of unboundedness is known.  Numerical evidence points in the direction of {bounded} Widom factors, at least in the case of Jordan regions.  But no proof is currently available.
\begin{openproblem}
Does there exist a compact set $\fre\subset\C$ with finitely many components for which the Widom factors are unbounded?
\end{openproblem}
The above mentioned results are qualitative in nature as the involved constants are large and their dependence on the set $\fre$ is rather implicit. Going in the other direction and assuming smoothness of the components of $\fre$ typically yields more explicit constants and even precise asymptotics for the Widom factors and the Chebyshev polynomials.

Suppose now that $\fre\subset\bbC$ is a finite disjoint union of $C^{2+}$ arcs and/or Jordan regions with $C^{2+}$ boundaries. Assume also that the weight function $w$ is supported on the boundary of $\fre$. Under these assumptions, Widom \cite[Sect.~11]{Wid69} obtained the asymptotic upper bound
\begin{align}
\label{wUB2}
  \limsup_{n\to\infty} \|w T_{n,w}\|_\fre / \ca(\fre)^n \le
2S(w)\exp\bigl[PW(\fre)\bigr],
\end{align}
compare with \eqref{TWub}.
This asymptotic bound is sharp within the class of real sets (i.e., $\fre$ consisting only of arcs lying on the real line).  However, in the case of $\fre$ consisting only of regions, Widom \cite[Sect.~8]{Wid69} established the improved asymptotic upper bound
\begin{align}
\label{wUB1}
  \limsup_{n\to\infty} \|w T_{n,w}\|_\fre / \ca(\fre)^n \le S(w)\exp\bigl[PW(\fre)\bigr].
\end{align}
More remarkably, in that case Widom showed that we have Szeg\H{o}--Widom asymptotics for the weighted Chebyshev polynomias $T_{n,w}$ and their norms $\|w T_{n,w}\|_\fre$ (i.e., the weighted analogs of \eqref{W-surmise} and \eqref{W-asympt}).
The improved asymptotic bound \eqref{wUB1} is also sharp; in fact, by \eqref{wLB}, equality is attained when $\fre$ consists of a single region since in that case the Green's function has no critical points and thus $PW(\fre)=0$.

For special subsets of the complex plane, we also have non-asymptotic upper bounds that parallel the real case. 
The following two results are taken from \cite{CSZ4}.
\begin{theorem}
If $\fre\subset\bbC$ is a solid lemniscate, that is,
\begin{equation}
  \fre=\{z\in\bbC \st |P(z)|\le\al\}
\end{equation}
for some polynomial $P$ of degree $k\ge1$ and $\al>0$, then 
\begin{equation}
  \|T_n\|_\fre \le K\cdot \ca(\fre)^n,
\end{equation}
where the constant $K$ is given by
\begin{equation}
  K = \max_{j=0,\dots,k-1} W_j(\fre).
\end{equation}
\end{theorem}

The other special case is motivated by an old result of Faber \cite{Fab1919} stating that the Chebyshev polynomials of an ellipse are the same as the ones for the interval between the two foci. This in particular leads to explicit values of the Widom factors for ellipses.  By further developing the results of Fischer \cite{Fis92} for two intervals, one can produce general results for level sets of the Green's function.
\begin{theorem}
If $\fre_0\subset\bbR$ is a PW set and
\begin{equation}
  \fre = \{z\in\bbC \st G(z)\le\alpha\}
\end{equation}
for some $\alpha > 0$, then 
\begin{equation}
  \bigl\|T^{(\fre)}_n\bigr\|_\fre \le \bigl(1+e^{-n\alpha}\bigr)\exp\bigl[PW(\fre_0)\bigr]\ca(\fre)^n.
\end{equation}
In addition, if $\fre_0$ is a period-$n$ set then 
the Chebyshev polynomials of degree $nk$ for the sets $\fre$ and $\fre_0$ coincide and
\begin{equation}
  \bigl\|T^{(\fre)}_{nk}\bigr\|_\fre = \cosh(nk\alpha)\bigl\|T^{(\fre_0)}_{nk}\bigr\|_{\fre_0}.
\end{equation}
\end{theorem}

\medskip

As mentioned above, in the case where $\fre$ consists of finitely many $C^{2+}$ Jordan regions, Widom obtained both Szeg\H{o}--Widom asymptotics and asymptotics of the Widom factors for the weighted Chebyshev polynomials. In the case of arcs, however, very little is known.  For weighted Chebyshev polynomials on finitely many interval (i.e., in the special case of arcs lying on the real line), Widom \cite[Sect.~11]{Wid69} merely established asymptotics of the Widom factors and conjectured the corresponding Szeg\H{o}--Widom asymptotics for the polynomials.  This conjecture was proven in the unweighted case in \cite{CSZ1}, but remains open for the weighted case.

In \cite{TY15}, Totik--Yuditskii extended the asymptotics of the Widom factors for weighted Chebyshev polynomials to the case of $\fre$ consisting of finitely many intervals and $C^{2+}$ Jordan regions symmetric with respect to the real line.  Yet, the case of sets consisting of finitely many smooth components some or all of which are arcs in general position in the complex plane has proven to be much more difficult.  Widom made conjectures regarding that case, but subsequent works \cite{TD91,TY15,Eic17} have shown that these conjectures are incorrect.  In particular, Widom expected that generically the asymptotic upper bound \eqref{wUB2} is attained for sets $\fre$ with finitely many smooth components when at least one of them is an arc.  While this was shown to be false in \cite{TY15}, the same work \cite{TY15} and in the unweighted case \cite{TD91} also showed that Widom was qualitatively correct in expecting larger asymptotics when an arc component is present.  In addition, for unweighted Chebyshev polynomials it was shown in \cite{Tot14} that sets $\fre$ containing an arc lead to an increased lower bound
\begin{align}
  \|T_{n}\|_\fre \ge (1+\be) \ca(\fre)^n, \quad n\ge1,
\end{align}
for some $\be>0$ that depends only on $\fre$ (cf. \eqref{Sze1}).

So far, the only nontrivial example of an arc for which the asymptotics is known is a single arc on the unit circle. In that case, Widom expected the asymptotics to be the same as for an interval. However, it was observed in \cite{TD91} that for the circular arc $\fre=\bigl\{e^{i\theta}:\theta\in[-\al,\al]\bigr\}$ (with $0<\al<\pi$), the unweighted Widom factors obey the asymptotics
\begin{equation}
\label{circ-arc-lim}
  \lim_{n\to\infty} W_n(\fre) = 1+\cos\Big(\frac{\al}2\Big).
\end{equation}
This shows that the case of a circular arc continuously interpolates between the case of a region (e.g., $W_n(D)\equiv1$ for a closed disk $D$) and the case of a flat arc (e.g., $W_n(I)\equiv2$ for an interval $I$).
In addition, it was shown in \cite{SZ21} that the Widom factors for a circular arc are strictly monotone increasing. The Szeg\H{o}--Widom asymptotics for the unweighted Chebyshev polynomials of a circular arc was derived by Eichinger \cite{Eic17} and the behavior is indeed different from the case of an interval.

At this point, we also mention a curious observation made in \cite[Thm.~5.1]{AZ20}.  For the polynomials orthogonal with respect to the equilibrium measure $d\rho_\fre$ on a circular arc $\fre$, the square of the associated $L^2$-Widom factors have the same asymptotics as $W_n(\fre)$ in \eqref{circ-arc-lim}.  This suggests that the two quantities might also coincide for other smooth arcs in the complex plane.  Since the asymptotics of the $L^2$-Widom factors for a $C^{2+}$ arc is known (see \cite{Alp22,Wid69}), we are led to the following conjecture:
\begin{conjecture}
If $\fre$ is a smooth arc in the complex plane, then
\begin{equation}
\label{arc-conj}
  \lim_{n\to\infty} W_n(\fre) = 2\pi S(w_\fre)\ca(\fre),
\end{equation}
where $w_\fre=\frac{1}{2\pi}\bigl(\frac{\pd G}{\pd n_+}+\frac{\pd G}{\pd n_-}\bigr)$ is the density of the equilibrium measure $d\rho_\fre$ with respect to arc-length.
\end{conjecture}
 For $C^{2+}$ arcs, Alpan \cite{Alp22} showed that the conjectured asymptotic value satisfies
\begin{equation}
1<2\pi S(w_\fre)\ca(\fre)\le2
\end{equation}
with the upper bound being strict if and only if $\frac{\pd G}{\pd n_+}(z) \neq \frac{\pd G}{\pd n_-}(z)$ for some non-endpoint $z\in\fre$.  The latter holds, for example, for non-analytic arcs.  Partial progress towards the above conjecture is also reported in \cite[Thm.~1.3]{Alp22} where the asymptotic upper bound \eqref{wUB2} is improved by replacing the constant $2$ with the smaller constant $2\sqrt{\pi S(w_\fre)\ca(\fre)}$.


The study of Chebyshev polynomials for subsets of the complex plane has another interesting and challenging direction which concerns the asymptotic behavior of their zeros.  Let $w_1,\dots,w_n$ be the zeros of $T_n$ counting multiplicity and denote by
\begin{equation}
 d\mu_n=\frac{1}{n}\sum_{j=1}^{n}\delta_{w_j}
\end{equation}
the normalized zero-counting measure for $T_n$.  The limit points of $\{d\mu_n\}_{n=1}^\infty$ as $n\to\infty$ are called \emph{density of Chebyshev zeros} for $\fre$.

In \cite{Wid67}, Widom proved 
that for any closed subset $S$ of $\Om$, the unbounded component of $(\C\cup\{\infty\})\setminus\fre$, there is an upper bound on the number of zeros of $T_n$ in $S$ which depends only on $S$ and not on $n$.  This implies the following general result on the density of Chebyshev zeros as stated in \cite{CSZ4}.

\begin{theorem}
Any limit point $d\mu_\infty$ of the zero-counting measures $d\mu_n$ is supported in the polynomial convex hull of $\fre$.  Moreover, for all $z\in\Om$ we have that
\begin{equation}
  \int \log|z-w|\, d\mu_\infty(w) = \int \log|z-w|\, d\rho_\fre(w).
\end{equation}
\end{theorem}

This theorem says that $d\rho_\fre$ is the balayage (see, e.g., \cite[Sect.~II.4]{SafTot97}) of $d\mu_\infty$ onto $\partial\fre$, equivalently, the balayage of $d\mu_n$ converges to $d\mu_\fre$;  ideas that go back at least to Mhaskar--Saff \cite{MhaSaf91}.  It is an intriguing question to understand whether or not the zero-counting measures $d\mu_n$ (or some subsequence thereof) converge to the equilibrium measure $d\rho_\fre$.  In \cite{ST90}, Saff--Totik proved the following result.

\begin{theorem} 
Let $\fre\subset\bbC$ be a compact set with connected interior and complement. Then:
\begin{enumerate}[$(a)$]
\item If $\fre$ is an analytic Jordan region (i.e., $\partial \fre$ is an analytic simple curve), then there is a neighborhood $U$ of $\partial\fre$ so that for all large $n$, $T_n$ has no zeros in $U$.

\item If $\partial\fre$ has a neighborhood $U$ and there is a sequence $n_j\to\infty$ so that $\mu_{n_j}(U)\to0$, then $\fre$ is an analytic Jordan region.
\end{enumerate}
\end{theorem}
Accordingly, for analytic Jordan regions the equilibrium measure is never a density of Chebyshev zeros and one may start wondering where these densities are supported.  Interestingly, and around the same time, Widom \cite{Wid90} had a similar result for nonselfadjoint Toeplitz matrices and Faber polynomials of the second kind.

In the complete other direction, Blatt--Saff--Simkani \cite{BSS88} proved the following result.
\begin{theorem}
Let $\fre\subset\bbC$ be a polynomially convex set with empty interior.  Then, as $n\to\infty$, the Chebyshev zero-counting measures $d\mu_n$ converge weakly to $d\rho_\fre$.
\end{theorem}

As explained below, there are also local versions of the above two theorems (see \cite{CSZ4} for proofs).
\begin{theorem}
Let $\fre\subset\bbC$ be a polynomially convex set and suppose $U\subset\bbC$ is an open connected set with connected complement so that $U\cap\partial\fre$ is a continuous arc that divides $U$ into two pieces, $\fre^{\inte}\cap U$ and $(\bbC\setminus\fre)\cap U$.  If $M_n(U)$ denotes the number of zeros of $T_n$ in $U$ and
\begin{equation}
  \liminf_{n\to\infty} \frac{M_n(U)}{n} = 0
\end{equation}
then
\begin{equation}
  U\cap\partial\fre \mbox{ is an analytic arc}.
\end{equation}
\end{theorem}
It readily follows that if $\fre$ is a Jordan region whose boundary curve is piecewise analytic but not analytic at some corner points, then at least these corner points are points of density for the zeros of $T_n$. Moreover, if $\partial\fre$ is nowhere analytic then all of the boundary points are points of density for the zeros.
In that light, it might be tempting to expect that the zero-counting measures $d\mu_n$ converge to the equilibrium measure $d\rho_\fre$ whenever $\partial\fre$ is nowhere analytic --- and this was conjectured in \cite{CSZ4}. However, Totik \cite{Tot22} recently disproved such a statement (which was also considered by Widom \cite{Wid94} in the context of nonselfadjoint Toeplitz matrices).

Nevertheless, local convergence to the equilibrium measure can be proved in some cases.

\begin{theorem}
Let $\fre\subset\bbC$ be a polynomially convex set and suppose $U\subset\bbC$ is an open connected set whose complement is also connected.  Assume that $\ca(U\cap\fre)>0$ but that $U\cap\fre$ has two-dimensional Lebesgue measure zero.  Then, as $n\to\infty$, the zero-counting measures $d\mu_n$ restricted to $U$ converge weakly to the equilibrium measure $d\rho_\fre$ restricted to $U$.
\end{theorem}

Another interesting result on convergence to the equilibrium measure is given by Saff--Stylianopoulos \cite{SS15}. They prove that if $\partial\fre$ has an inward pointing corner (more generally, a non-convex type singularity), then the zero-counting measures $d\mu_n$ always converge weakly to $d\rho_\fre$.  For example, if $\fre$ is a non-convex polygon then their hypothesis holds.  The case of convex polygons, on the other hand, leads to an interesting open problem.
\begin{openproblem}
What are the density of Chebyshev zeros when $\fre$ is a convex polygon?
\end{openproblem}
This is not even known for the equilateral triangle, although numerical computations present some evidence for convergence 
to the equilibrium measure. The other natural candidate for the limit points of zeros is the skeleton consisting of the line segments from the centroid of the triangle to the vertices.



\end{document}